# THE LINKS-GOULD INVARIANT AS A CLASSICAL GENERALIZATION OF THE ALEXANDER POLYNOMIAL ?

BEN-MICHAEL KOHLI

Abstract. In this paper we conjecture that the Links-Gould invariant of links, that we know is a generalization of the Alexander-Conway polynomial, shares some of its classical features. In particular it seems to give a lower bound for the genus of links and to provide a criterion for fiberedness of knots. We give some evidence for these two assumptions.

## Contents

## Introduction

The Links-Gould invariant $LG(L; t_0, t_1)$ is a two variable polynomial quantum invariant. It is derived from a one parameter family of representations of quantum superalgebra $U_q(gl(2|1))$ [5, 18]. It is part of a larger family of Links-Gould invariants $LG^{m,n}$, $m, n \in \mathbb{N}^*$ [7].

It is worth noticing that the Alexander-Conway polynomial of a link $\Delta_L$ can be recovered from $LG$ in at least two ways. David de Wit, Atsushi Ishii and Jon Links showed [4]

$$LG(L; t_0, -t_0^{-1}) = \Delta_L(t_0^2).$$

The square of the Alexander polynomial can also be obtained evaluating $LG$ [16, 17]

$$LG(L; t_0, t_0^{-1}) = \Delta_L(t_0)^2.$$

Knowing this, it is natural to wonder :

**Question 0.1.** *Are there properties of $\Delta$ that extend to $LG$ ?*

In particular, if $\Delta$ can be seen as a quantum invariant [23, 25, 30], it is in essence a classical invariant derived from a presentation matrix of the first homology group of the infinite cyclic covering of the complement of a given link in $S^3$ [1]. Therefore we can ask if some of $\Delta$'s homological properties extend to $LG$.

In that spirit, let us recall that a conjecture by Ishii [12] states :

**Conjecture 0.2.** *The LG polynomial $LG(K; t_0, t_1) = \sum_{i,j} a_{ij} t_0^i t_1^j$ of an alternating knot $K$ is "alternating", that is : $a_{ij}a_{kl} \geqslant 0$ if $i + j + k + l$ is even, and $a_{ij}a_{kl} \leqslant 0$ otherwise.*

2010 *Mathematics Subject Classification.* 57M27.

Though if you think about it this is not a straightforward generalization of the well known similar property for the Alexander polynomial [20, 3] for either of the two evaluations we know $t_0t_1 = 1$ and $t_0t_1 = -1$, it still can be thought of as the trace of a similar behavior.

In the following we give evidence for more positive answers to Question 0.1. We conjecture that the span of the $LG$ invariant is a lower bound for the genus of a link.

**Conjecture 0.3.** *Set $L$ a link in $S^3$ and $\mu$ the number of its components.*

- **I:** *$span(LG(L; t_0, t_1)) \leqslant 2(2g(L) + \mu - 1)$,*
- **II:** *If $L$ is alternating, then inequality* **I** *is an equality.*

We also conjecture that for fibered knots, there are conditions on the leading coefficients of the $LG$ polynomial.

**Conjecture 0.4.** *Set $K$ a knot in $S^3$.*

- **I:** *If $K$ is fibered then $LG(K)$ is monic,*
- **II:** *If $K$ is alternating, the converse is true as well.*

We base these conjectures on computations for the first prime knots and on partial skein relations for $LG$ that allow its evaluation on various infinite families of links. Notice that if the genus conjecture were true, $LG$ would systematically give a better lower bound for the genus of a link than the one given by the Alexander invariant. Also, the criterion we conjecture for fibered knots would refine the well known similar statement for $\Delta$.

A proof of these two statements would show *quantum* invariant $LG$ can be used to find information on the geometry of links.

## 1. The Links-Gould invariant and the genus of links

1.1. **Definitions.** We recall the definition of the Links-Gould invariant of links. We then give the definition of the genus of a link and recall the connection between the genus of a link and its Alexander invariant.

**Definition 1.1.** Set $\mathbb{K} := \mathbb{C}(t_0^{\pm\frac{1}{2}}, t_1^{\pm\frac{1}{2}})$. Let $W =< e_1, \ldots, e_4 >$ be a four-dimensional $\mathbb{K}$-vector space. The following linear map $R$, expressed in basis $(e_1 \otimes e_1, e_1 \otimes e_2, e_1 \otimes e_3, e_1 \otimes e_4, e_2 \otimes e_1, e_2 \otimes e_2, e_2 \otimes e_3, \ldots)$, is an automorphism of $W \otimes W$ and an R-matrix ([5], p.186) :

$$\begin{pmatrix}
t_0 & . & . & . & . & . & . & . & . & . & . & . & . & . & . & . \\
. & . & . & . & t_0^{1/2} & . & . & . & . & . & . & . & . & . & . & . \\
. & . & . & . & . & . & . & . & t_0^{1/2} & . & . & . & . & . & . & . \\
. & . & . & . & . & . & . & . & . & . & . & . & 1 & . & . & . \\
. & t_0^{1/2} & . & . & t_0-1 & . & . & . & . & . & . & . & . & . & . & . \\
. & . & . & . & . & -1 & . & . & . & . & . & . & . & . & . & . \\
. & . & . & . & . & . & t_0t_1-1 & . & . & -t_0^{1/2}t_1^{1/2} & . & . & -t_0^{1/2}t_1^{1/2}Y & . & . & . \\
. & . & . & . & . & . & . & . & . & . & . & . & . & t_1^{1/2} & . & . \\
. & . & t_0^{1/2} & . & . & . & . & . & t_0-1 & . & . & . & . & . & . & . \\
. & . & . & . & . & . & -t_0^{1/2}t_1^{1/2} & . & . & . & . & . & Y & . & . & . \\
. & . & . & . & . & . & . & . & . & . & -1 & . & . & . & . & . \\
. & . & . & . & . & . & . & . & . & . & . & . & . & . & t_1^{1/2} & . \\
. & . & . & 1 & . & . & -t_0^{1/2}t_1^{1/2}Y & . & . & Y & . & . & Y^2 & . & . & . \\
. & . & . & . & . & . & . & t_1^{1/2} & . & . & . & . & . & t_1-1 & . & . \\
. & . & . & . & . & . & . & . & . & . & . & t_1^{1/2} & . & . & t_1-1 & . \\
. & . & . & . & . & . & . & . & . & . & . & . & . & . & . & t_1
\end{pmatrix}$$

where $Y = ((t_0 - 1)(1 - t_1))^{1/2}$.

We denote by $b_R^n$ the representation of braid group $B_n$ derived from this R-matrix. It is given by the standard formula

$$b_R^n(\sigma_i) = id_W^{\otimes i-1} \otimes R \otimes id_W^{\otimes n-i-1} \;,\; i = 1, \ldots, n-1.$$

*Remark* 1.2. Compared to the R-matrix used in [16], $R$ is multiplied by $-1$. We chose this convention here for it is the one Ishii uses in [12, 13], two papers that are interesting for our study.

**Theorem 1.3.** *Let $L$ be an oriented link, and $b \in B_n$ a braid with closure $L$. Set $\mu$ the linear map defined by*

$$\mu = \begin{pmatrix} t_0^{-1} & . & . & . \\ . & -t_1 & . & . \\ . & . & -t_0^{-1} & . \\ . & . & . & t_1 \end{pmatrix} \in End(W).$$

*Then :*
*1) there exists an element $c \in \mathbb{K}$ such that $trace_{2,3,\ldots,n}((id_W \otimes \mu^{\otimes n-1}) \circ b_R^n(b)) = c.id_W$,*
*2) $c$ is an oriented link invariant called Links-Gould invariant of $L$. We will denote it by $LG(L; t_0, t_1)$.*

*Remark* 1.4. In fact $LG(L; t_0, t_1) \in \mathbb{Z}[t_0^{\pm 1}, t_1^{\pm 1}]$ [12].

*Remark* 1.5. With the notations we use, $LG(L; q^{-2\alpha}, q^{2\alpha+2})$ is the Links-Gould invariant introduced in [5], using a one parameter family of representations of quantum superalgebra $U_q(gl(2|1))$.

**Definition 1.6.** (Seifert surface for a link) Set $L$ a link in $S^3$. A Seifert surface for $L$ is a compact, connected, orientable surface $\Sigma \subset S^3$ such that $\partial\Sigma = L$.

Such a surface exists for any link according to Seifert's algorithm [26].

*Remark* 1.7. Any Seifert surface $\Sigma$ being connected and orientable, one can define the genus $g(\Sigma)$ of $\Sigma$ :

$$\chi(\Sigma) = 2 - 2g(\Sigma) - \mu$$

where $\chi(\Sigma)$ is the Euler characteristic of $\Sigma$ and $\mu$ is the number of components of link $L$.

**Definition 1.8.** (genus of a link) Let $L$ be a link in $S^3$. The genus $g(L)$ of $L$ is

$$g(L) := min\{g(\Sigma),\ \Sigma \text{ Seifert surface for } L\}.$$

**Definition 1.9.** (Alexander polynomial of a link) Set $L$ a link in $S^3$ with $\mu$ components and choose $\Sigma$ a Seifert surface for $L$. Then $H_1(\Sigma, \mathbb{Z})$ is a free abelian group of rank $1 - \chi(\Sigma) = 2g(\Sigma) + \mu - 1$. If $v_{ij}$ is the linking number in $S^3$ of the $i^{th}$ generator of $H_1(\Sigma, \mathbb{Z})$ with the pushoff of the $j^{th}$ generator, then $V = (v_{ij})$ is a Seifert matrix for $L$. The Alexander polynomial is computed from such a Seifert matrix setting

$$\Delta_L(t) = det(tV - {}^tV) \in \mathbb{Z}[t, t^{-1}].$$

With this definition, $\Delta_L$ is determined up to multiplication by $\pm t^n$, $n \in \mathbb{Z}$. The standard Alexander normalization consists in picking the representative with positive constant term. The Alexander-Conway normalization corresponds, at least in the case of a knot $K$, to choosing the symmetric Laurent polynomial with $\Delta_K(1) = 1$.

**Proposition 1.10.** *For any link $L$, $deg(\Delta_L(t)) \leqslant 2g(L) + \mu - 1$.*

So the degree of $\Delta_L$ gives a lower bound on the genus of link $L$.

### 1.2. **The genus conjecture.**

We believe that Proposition 1.10 can be extended to the similar statement expressed in Conjecture 0.3. We will explain how and why it would be an extension of 1.10. The goal of section 2 is to give a range of evidence to support that conjecture.

**Definition 1.11.** Set $P \in \mathbb{Z}[t_0^{\pm 1}, t_1^{\pm 1}]$. For $(n, m) \in \mathbb{Z}^2$, we define

$$deg(t_0^n t_1^m) := n - m.$$

For a general $P = \sum_{i,j \in \mathbb{Z}} a_{ij} t_0^i t_1^j$, we can extend that definition to introduce the span of $P$ :

$$span(P) := max\{deg(t_0^i t_1^j),\ (i,j) \in \mathbb{Z}^2 \mid a_{ij} \neq 0\} - min\{deg(t_0^i t_1^j),\ (i,j) \in \mathbb{Z}^2 \mid a_{ij} \neq 0\}.$$

*Remark* 1.12. The span satisfies the usual elementary degree properties :

$$span(PQ) = span(P) + span(Q),$$

$span(P+Q) \leqslant max\{span(P), span(Q)\}$ **if** P and Q are symmetric Laurent polynomials.

Conjecture 0.3 generalizes Proposition 1.10 since this well known result shows Conjecture 0.3 is true when $t_0 t_1 = 1$ and $t_0 t_1 = -1$ via the evaluations we already mentionned

$$LG(L; t_0, -t_0^{-1}) = \Delta_L(t_0^2)\ ;\ LG(L; t_0, t_0^{-1}) = \Delta_L(t_0)^2.$$

These evaluations also explain why our definition for the span was natural to try and push the lower bound a little further.

**Proposition 1.13.**

**1 :** *In Conjecture 0.3,* **I** *implies* **II***,*

**2 :** *If Conjecture 0.3 is true, it systematically improves the lower bound for the genus provided by* $\Delta$ *:*

$$\textit{for any } L \textit{ link, } 2deg(\Delta_L(t)) \leqslant span(LG(L; t_0, t_1)).$$

*Moreover, there are links where*

$$2deg(\Delta_L(t)) < span(LG(L; t_0, t_1)).$$

*Proof.* Since $LG(L; t_0, -t_0^{-1}) = \Delta_L(t_0^2)$, $2deg(\Delta_L(t)) = deg(\Delta_L(t^2)) = span(LG(L; t, -t^{-1}))$. So to prove **2**, we wish to show

$$span(LG(L; t, -t^{-1})) \leqslant span(LG(L; t_0, t_1)).$$

If we denote $LG(L; t_0, t_1) = \sum_{\substack{i,j \in \mathbb{Z}/ \\ a_{ij} \neq 0}} a_{ij} t_0^i t_1^j$, then

$$LG(L; t, -t^{-1}) = \sum_{k \in \mathbb{Z}} \Big( \sum_{\substack{i,j \in \mathbb{Z}/ \\ a_{ij} \neq 0 \\ i-j=k}} a_{ij} (-1)^j \Big) t^k.$$

This clearly shows that if the coefficient in front of $t^k$ in $LG(L; t, -t^{-1})$ is non zero, then there is at least one non zero coefficient in front of a monomial of degree $k$ in the expression of $LG(L; t_0, t_1)$, which yields **2**. Moreover, some examples where the equality does not hold are given in Proposition 2.1.

Now suppose **I** holds for any link and set $L$ an alternating link. Then [3], Theorem 3.5, states

$$deg(\Delta_L(t)) = 2g(L) + \mu - 1.$$

So we have the following inequality chain :

$$\begin{aligned} 2.(2g(L) + \mu - 1) &\geqslant span(LG(L; t_0, t_1)) && \text{(Conjecture 0.3)} \\ &\geqslant 2deg(\Delta_L(t)) && \text{(point } \mathbf{2}\text{)} \\ &= 2.(2g(L) + \mu - 1) && \text{(reference [3])} \end{aligned}$$

CQFD

## 2. Evidence supporting the genus conjecture

We wish to give evidence of the likeliness of Conjecture 0.3. In particular we verify the bound for small prime knots, prove it for several infinite families of knots and links and verify that the genus conjecture holds on an untwisted Whitehead double of the trefoil knot, which is a counter example due to Hugh Morton in a ressembling situation we will explain.

### 2.1. **Less than 13 crossing prime knots.**

When we consider knots, $\mu$ is equal to 1 and the inequality becomes

$$span(LG(K;t_0,t_1)) \leqslant 4g(K).$$

We tested that inequality on all prime knots with less than 12 crossings, and on a large selection of non alternating prime knots with 13 crossings. To do that, we used the computations of $LG$ for prime knots one can access via David de Wit's LINKS-GOULD EXPLORER [8]. To find genus information up to 12 crossings, we used Cha and Livingston's KNOT-INFO [2]. For non alternating 13 crossing prime knots, data is obtained from Stoimenow's website KNOT DATA TABLES [29]. Knots are listed with respect to the HTW ordering for tables of prime knots of up to 16 crossings [11].

The reason why we did not test all non alternating 13 crossing prime knots is explained in [6] : The LINKS-GOULD EXPLORER's database contains evaluations only for $LG$ of knots with string index at most 5, and from time to time 6 or 7. Indeed, the memory required increases dramatically with braid width. This still provides values for $LG$ for 2096 non alternating prime knots with 13 crossings among the 5110 which exist.

**Proposition 2.1.**

1 : *Conjecture 0.3 holds for every knot tested. In particular, for all alternating knots tested, the equality holds.*

2 : *For all prime knots with less than 10 crossings, the span of LG exactly is 4 times the genus of the knot.*

3 : *The list of prime knots with 11 or 12 crossings where there is no equality is the following :* $11^N_{34}$, $11^N_{42}$, $11^N_{45}$, $11^N_{67}$, $11^N_{73}$, $11^N_{97}$, $11^N_{152}$, $12^N_{28}$, $12^N_{31}$, $12^N_{51}$, $12^N_{56}$, $12^N_{63}$, $12^N_{87}$, $12^N_{129}$, $12^N_{132}$, $12^N_{221}$, $12^N_{231}$, $12^N_{256}$, $12^N_{257}$, $12^N_{264}$, $12^N_{267}$, $12^N_{268}$, $12^N_{313}$, $12^N_{430}$, $12^N_{665}$, $12^N_{808}$, $12^N_{812}$.

4 : *Among these knots we get a more precise lower bound with LG than with $\Delta$ for several of them :* $11^N_{34}$, $11^N_{42}$, $11^N_{67}$, $11^N_{97}$, $12^N_{31}$, $12^N_{51}$, $12^N_{129}$, $12^N_{256}$, $12^N_{257}$, $12^N_{264}$, $12^N_{267}$, $12^N_{268}$, $12^N_{313}$, $12^N_{430}$, $12^N_{665}$, $12^N_{812}$.

5 : *Some of the spans are "half integers" in the sense that they are multiples of 2, as they necessarily are, but not multiples of 4 :* $11^N_{34}$, $11^N_{42}$, $11^N_{67}$, $11^N_{97}$, $12^N_{51}$, $12^N_{256}$, $12^N_{257}$, $12^N_{264}$, $12^N_{267}$, $12^N_{268}$, $12^N_{313}$, $12^N_{430}$, $12^N_{665}$, $12^N_{812}$.

*Remark* 2.2. Once it is proved to be true, the span inequality is meaningful when $2deg(\Delta_L(t)) < span(LG(L;t_0,t_1))$. However, as long as it remains a conjecture, the hard case for the inequality is when $deg(\Delta_L(t)) = 2g(L)+\mu-1$ precisely because there is no choice on the value of the span of $LG$ for it not to be a counter-example.

### 2.2. **The untwisted Whitehead double of the trefoil knot.**

Here we compute the Links-Gould polynomial of the untwisted double of the trefoil knot that is drawn is Figure 1. This is an interesting knot to study since it is a counterexample to a genus type bound for another generalization of the Alexander-Conway polynomial : the HOMFLY-PT polynomial. Precisely, Hugh Morton shows in [19], theorem 2, that the monomial of highest degree with respect to the Alexander variable in the HOMFLY-PT invariant gives a lower bound for $2\tilde{g}(L)+\mu-1$, where $\tilde{g}(L)$ is the *canonical genus* of link $L$. The double of the trefoil is the example Morton gives to show that degree is *not* in general a lower bound for $2g(L)+\mu-1$. There is no such problem here :

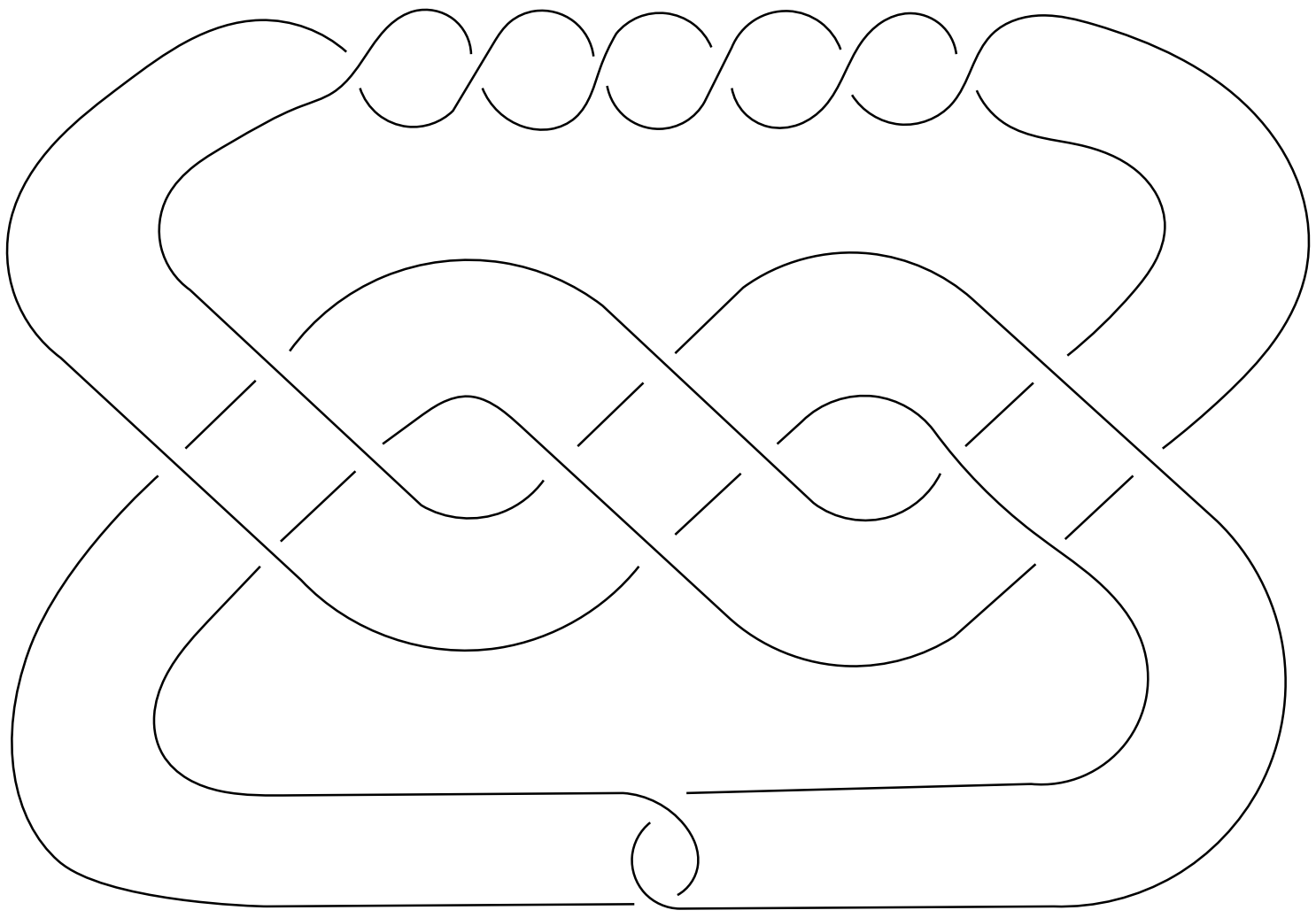

FIGURE 1. The untwisted Whitehead double of the right handed trefoil with a positive clasp.

**Proposition 2.3.** *The untwisted Whitehead double of the right handed trefoil with positive clasp $K_0$ has the following properties :*

- $g(K_0) = 1$,
- $\Delta_{K_0}(t) = 1$,
- *A braid presentation of $K_0$ in braid group $B_6$ as a word in the standard Artin generators $\sigma_1, \ldots, \sigma_5$ is :*
$$\sigma_4\sigma_3\sigma_2\sigma_3^{-1}\sigma_4^{-1}\sigma_5\sigma_3\sigma_2\sigma_1\sigma_2^{-1}\sigma_3^{-1}\sigma_5\sigma_4\sigma_3\sigma_4^{-1}\sigma_5^{-1}\sigma_3\sigma_2\sigma_3^{-1}\sigma_2\sigma_1^2\sigma_2^{-1},$$
- $LG(K_0; t_0, t_1) = 3 - 4t_1 + 2t_1^2 - 4t_0 + 6t_1t_0 - 2t_1^2t_0 + 2t_0^2 - 2t_1t_0^2 - 2t_1^2t_0^2 + 4t_1^3t_0^2 - 2t_1^4t_0^2 + 4t_1^2t_0^3 - 10t_1^3t_0^3 + 8t_1^4t_0^3 - 2t_1^5t_0^3 - 2t_1^2t_0^4 + 8t_1^3t_0^4 - 8t_1^4t_0^4 + 2t_1^5t_0^4 - 2t_1^3t_0^5 + 2t_1^4t_0^5 + 4t_1^5t_0^5 - 6t_1^6t_0^5 + 2t_1^7t_0^5 - 6t_1^5t_0^6 + 8t_1^6t_0^6 - 2t_1^7t_0^6 + 2t_1^5t_0^7 - 2t_1^6t_0^7.$

*So the span of $LG(K_0; t_0, t_1)$ is 4 and the span inequality is verified in this case.*

*Remark* 2.4. The value of $LG(K_0; t_0, t_1)$ was obtained by direct computation of the formula given by Theorem 1.3 with the R-matrix in Definition 1.1 using MATHEMATICA 10.

2.3. **Infinite families of knots where Conjecture 0.3 can be verified thanks to partial skein relations we know.** Here we verify the genus bound on several infinite families of knots or links. To do that, we will use basic Alexander-type properties of $LG$ we will recall, and partial skein relations that will make the computations practicable.

2.3.1. *Some properties of LG and useful skein relations.* To compute $LG$ for infinite families of knots, we need to have a more efficient way to evaluate it than simply using the formula in 1.3. We first recall some general facts about the $LG$ polynomial.

**Proposition 2.5.** *The Links-Gould polynomial satisfies the following properties :*

- $LG(\bigcirc) = 1$,
- *Denoting $L^*$ the reflexion of $L$, $LG(L^*; t_0, t_1) = LG(L; t_0^{-1}, t_1^{-1})$,*
- *We have the following symmetry : $LG(L; t_0, t_1) = LG(L; t_1, t_0)$. Indeed $LG$ does not detect inversion,*
- *For $L$ and $L'$ two links, denoting $L\#L'$ their connected sum : $LG(L\#L') = LG(L)LG(L')$,*
- *If $L = L' \sqcup L''$ is the split union of $L'$ and $L''$, then $LG(L) = 0$.*

*Proof.* For proofs of these facts, we refer the reader to [13, 5, 7]. CQFD

*Remark* 2.6. The last two points show that $LG$ and $\Delta$ behave similarily concerning sums and disjoint unions.

Let us also cite a list of skein relations that are known to be true for $LG$. Whether the associated skein module is generated by the unknot or not is a problem pointed out by Ishii [13]. It is to the best of our knowledge an open question.

**Proposition 2.7.** *LG verifies the following skein relations.*
*Skein relation* (1) :

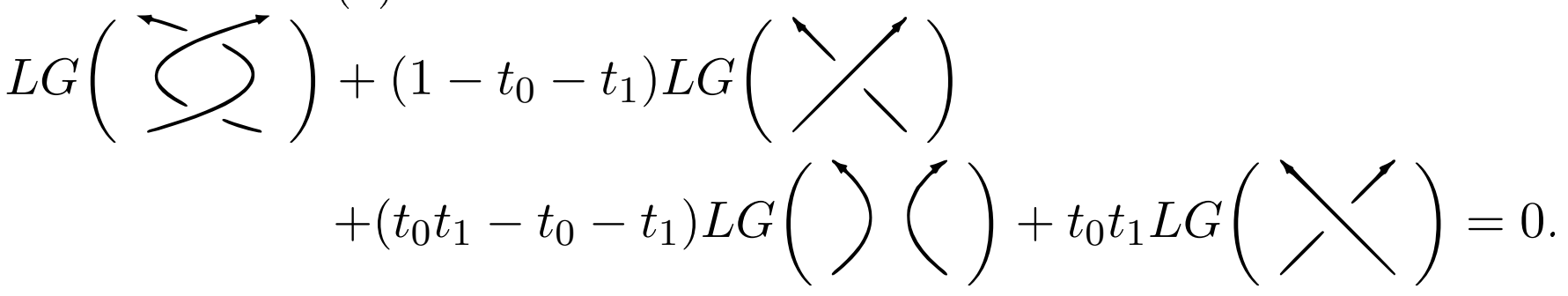

*Skein relation* (2) :

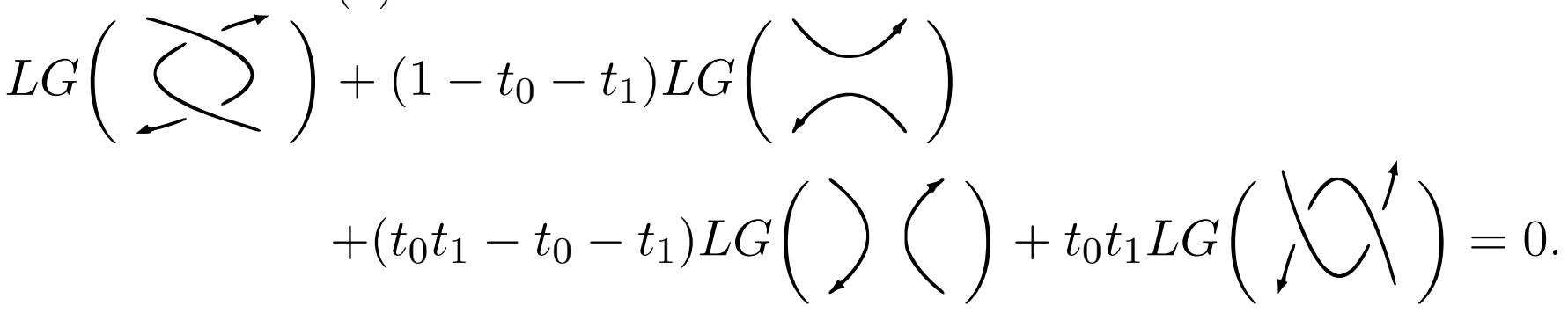

*Skein relation* (3) :

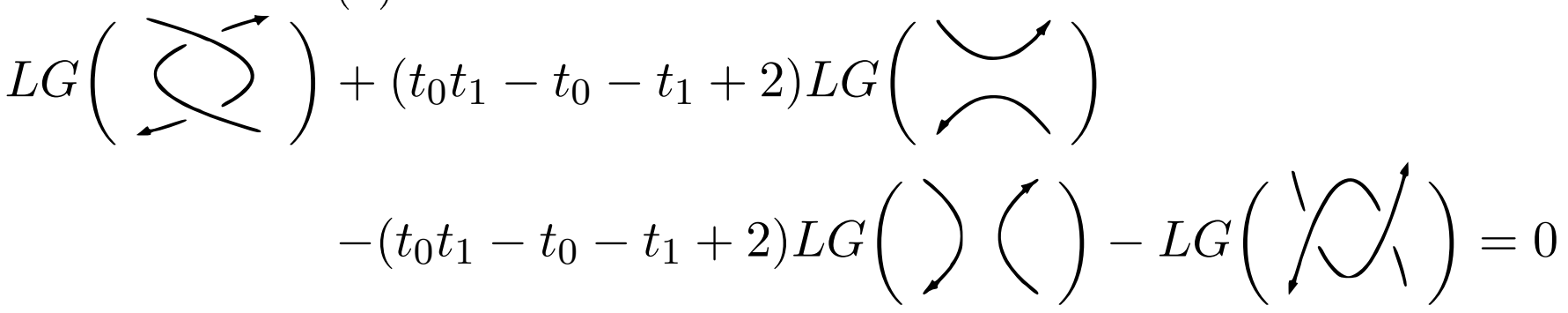

.

*Skein relation* (4) :

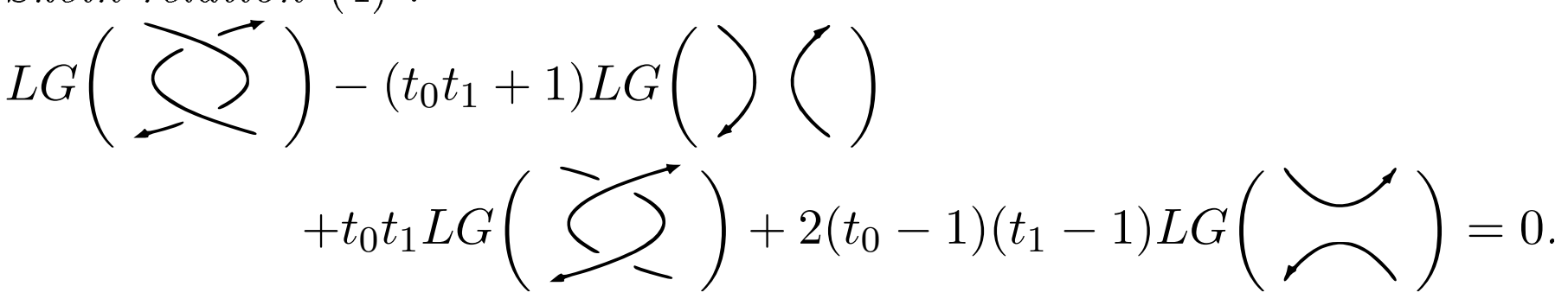

*Proof.* See [13, 5, 18]. CQFD

*Remark* 2.8. (1) and (2) are equivalent, by adding each time a well chosen tangle from the left.

*Remark* 2.9. As explained in [13], (4) is a consequence of (2) and (3).

*Remark* 2.10. Set $V$ the 4-dimensional irreducible $U_q(gl(2|1))$-module that gives rise to the Links-Gould invariant. Then the tensor product of two copies of $V$ decomposes with respect to the $U_q(gl(2|1))$-module structure.

$$V \otimes V = V_1 \oplus V_2 \oplus W, \text{ with } dimV_1,\ dimV_2 = 4 \text{ and } dimW = 8.$$

Morevover, $V_1$, $V_2$ and $W$ are non isomorphic irreducible $U_q(gl(2|1))$-modules. For details, see [10, 5]. Using this and denoting $A := U_q(gl(2|1))$, we have the following identification :

$$End_A(V \otimes V) \simeq End_A(V_1) \oplus End_A(V_2) \oplus End_A(W) \simeq \mathbb{C}(t_0^{\pm 1}, t_1^{\pm 1})^3.$$

In particular, for *any* three $(2,2)$-tangles such that the associated maps in $End_A(V \otimes V)$ are linearly independent, *any* other can be expressed as a linear combination of the first three. This potentially generates a great variety of skein relations for $LG$.

*Remark* 2.11. Using points 2 and 3 of Proposition 2.5, we can modify the previous skein relations : orientation of the strands, signs of the crossings. We will use these modified relations, though we will not write them down here.

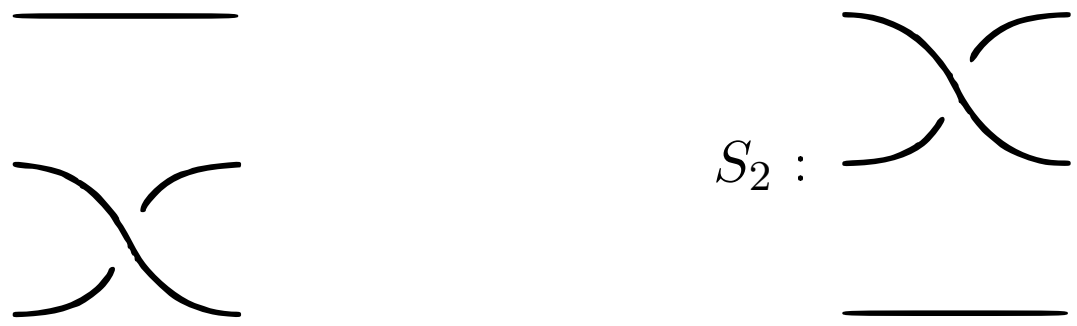

FIGURE 2. Generators $S_1$ and $S_2$ of braid group $B_3$.

**Corollary 2.12.** *Using notations in* [13], *LG satisfies the following skein relations :*

- $LG\Big($ 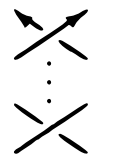 $\Big\}$ *n half twists* $\Big)$

$= \left(\frac{(-1)^n}{(t_0+1)(t_1+1)} + \frac{t_0^n}{(t_0+1)(t_0-t_1)} + \frac{t_1^n}{(t_1+1)(t_1-t_0)}\right) LG\Big($ 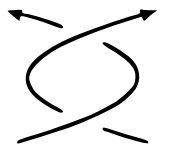 $\Big)$

$- \left(\frac{(-1)^n(t_0+t_1)}{(t_0+1)(t_1+1)} + \frac{t_0^n(t_1-1)}{(t_0+1)(t_0-t_1)} + \frac{t_1^n(t_0-1)}{(t_1+1)(t_1-t_0)}\right) LG\Big($ 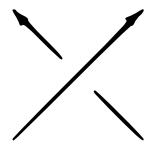 $\Big)$

$+ \left(\frac{(-1)^n t_0 t_1}{(t_0+1)(t_1+1)} - \frac{t_0^n t_1}{(t_0+1)(t_0-t_1)} - \frac{t_1^n t_0}{(t_1+1)(t_1-t_0)}\right) LG\Big($ 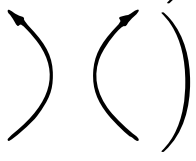 $\Big)$.

- $LG\Big($ 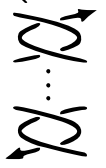 $\Big\}$ *n full twists* $\Big)$

$= \frac{t_0^n t_1^n - 1}{t_0 t_1 - 1} LG\Big($ 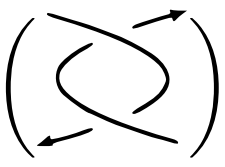 $\Big) + \left(1 - \frac{t_0^n t_1^n - 1}{t_0 t_1 - 1}\right) LG\Big($ 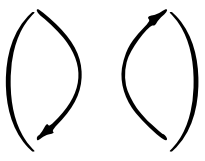 $\Big)$

$+ \frac{2(t_0-1)(t_1-1)}{t_0 t_1 - 1}\left(n - \frac{t_0^n t_1^n - 1}{t_0 t_1 - 1}\right) LG\Big($ 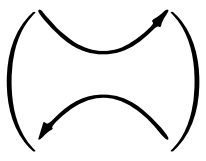 $\Big)$.

- $LG\Big($ 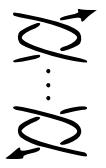 $\Big\}$ *n full twists* $\Big)$

$= a_1(n) LG\Big($ 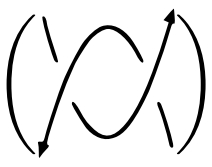 $\Big) + a_2(n) LG\Big($ 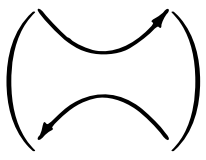 $\Big) + a_3(n) LG\Big($ 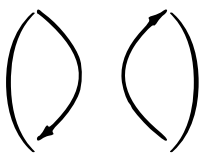 $\Big)$,

*where :*

$$a_1(n) = \frac{t_0^n t_1^n - 1}{t_0 t_1 - 1},$$
$$a_2(n) = \frac{2n(t_0-1)(t_1-1)}{t_0 t_1 - 1} - a_1(n)\left(\frac{(t_0 t_1+1)(t_0-1)(t_1-1)}{t_0 t_1 - 1} + 1\right),$$
$$a_3(n) = (t_0 - 1)(t_1 - 1)a_1(n) + 1.$$

We will now use all these properties to compute $LG$, or at least its span, on some infinite families of links.

2.3.2. *2-bridge links.* A 2-bridge link is a link with bridge number 2. As explained in [14], an oriented 2-bridge link can always be written in terms of the two generators $S_1$ and $S_2$ of 3-string braid group $B_3$. We use the notations one can find in [14, 13], setting $D(b_1, b_2, \ldots, b_m)$ the oriented 2-bridge link drawn in Figures 3 and 4. We can suppose $b_1, \ldots, b_m > 0$, thereby choosing an alternating diagram to represent any 2-bridge link. In particular, any 2-bridge link is alternating, so the inequality in Conjecture 0.3 should be an equality.

*Remark* 2.13. If $m$ is even, $D(b_1, \ldots, b_m)$ is a knot. It is a link with two components when $m$ is odd.

**Proposition 2.14.** *If* $b_i \neq 0$ *for any* $i$, *then*

$$g(D(b_1, \ldots, b_m)) = \frac{m - \mu + 1}{2}$$

*where* $\mu$ *is the number of components* [27].

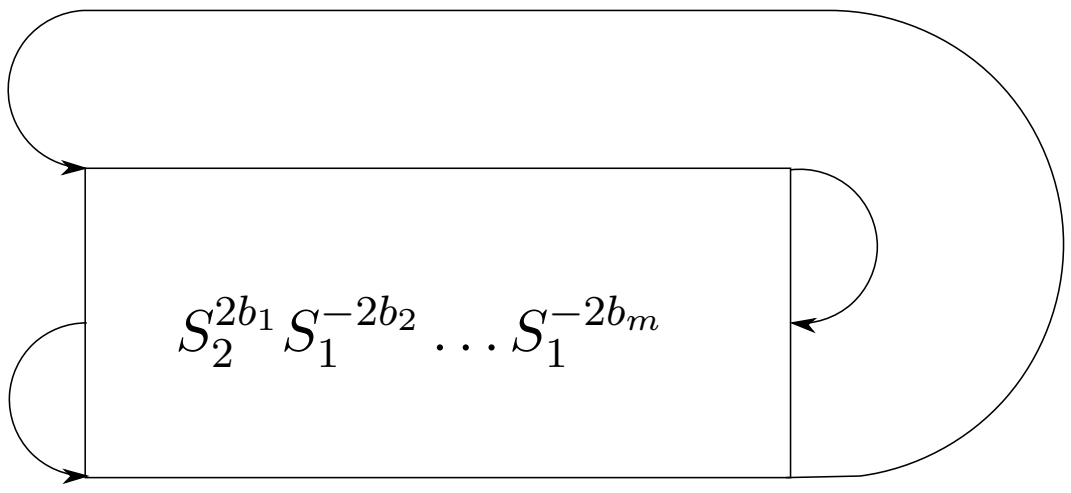

FIGURE 3. $D(b_1,\dots,b_m)$ when $m$ is even.

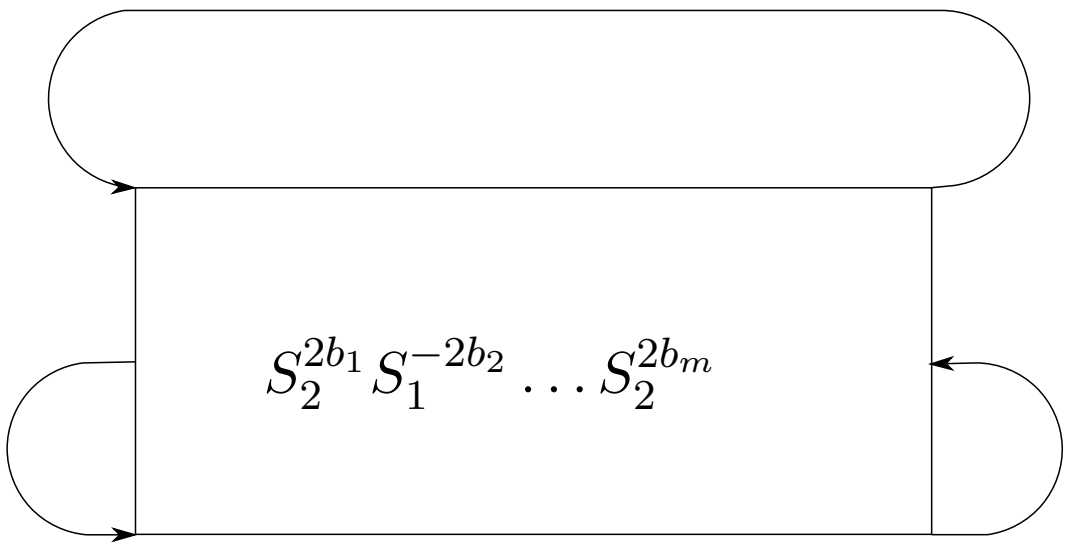

FIGURE 4. $D(b_1,\dots,b_m)$ when $m$ is odd.

Therefore, Conjecture 0.3 can be rephrased

$$span(LG(D(b_1,\dots,b_m))) = 2(2g(D(b_1,\dots,b_m)) + \mu - 1) = 2m.$$

**Theorem 2.15.** *For any* $b_1, b_2, \dots, b_m > 0$,

$$span(LG(D(b_1,\dots,b_m))) = 2m.$$

*Proof.* First we note that $a_1(n), a_2(n), a_3(n)$ are symmetric *polynomials* with respect to variables $t_0$ and $t_1$. We can compute the span in each case.

$$span(a_1(n)) = 0,\ span(a_2(n)) = 2,\ span(a_3(n)) = 2.$$

We will indicate by $\tilde{a_i}(n)$ the quantity $a_i(n)(t_0^{-1}, t_1^{-1})$. Let's prove the span equality by induction on $m$.

$\boxed{m=1}$ Using the mirror of skein relation 3 of Corollary 2.12 :

$$LG(D(b_1)) = \tilde{a_1}(b_1)LG(\bigcirc) + \tilde{a_2}(b_1)LG(\bigcirc) + \tilde{a_3}(b_1)LG(\bigcirc\bigcirc) = \tilde{a_1}(b_1) + \tilde{a_2}(b_1).$$

So the span of $LG(D(b_1))$ is 2.

$\boxed{m=2}$ Still using the same skein relation, we can compute $LG(D(b_1, b_2))$.

$$LG(D(b_1,b_2)) = a_1(b_2)LG(D(b_1-1)) + a_2(b_2)LG(D(b_1)) + a_3(b_2)LG(\bigcirc).$$

The second part of the sum has span $2+2=4$. The third part has span $2+0=2$. More care has to be taken with the first term, and in particular with $LG(D(b_1-1))$. If $b_1 - 1 > 0$, $LG(D(b_1-1))$ has span 2. In the other case, $D(0) = \bigcirc\bigcirc$ so $LG(D(0)) = 0$. So in any case the span of the sum is 4.

Let us now set $m \geqslant 3$ and suppose the equality stands for any $D(b_1,\dots,b_k)$ with $k \leqslant m-1$. For $D(b_1,\dots,b_m)$ we can apply skein relation 3 of 2.12 or its mirror image to the crossings that correspond to $S_1^{-2b_m}$ or $S_2^{2b_m}$ depending on whether $m$ is odd or even.

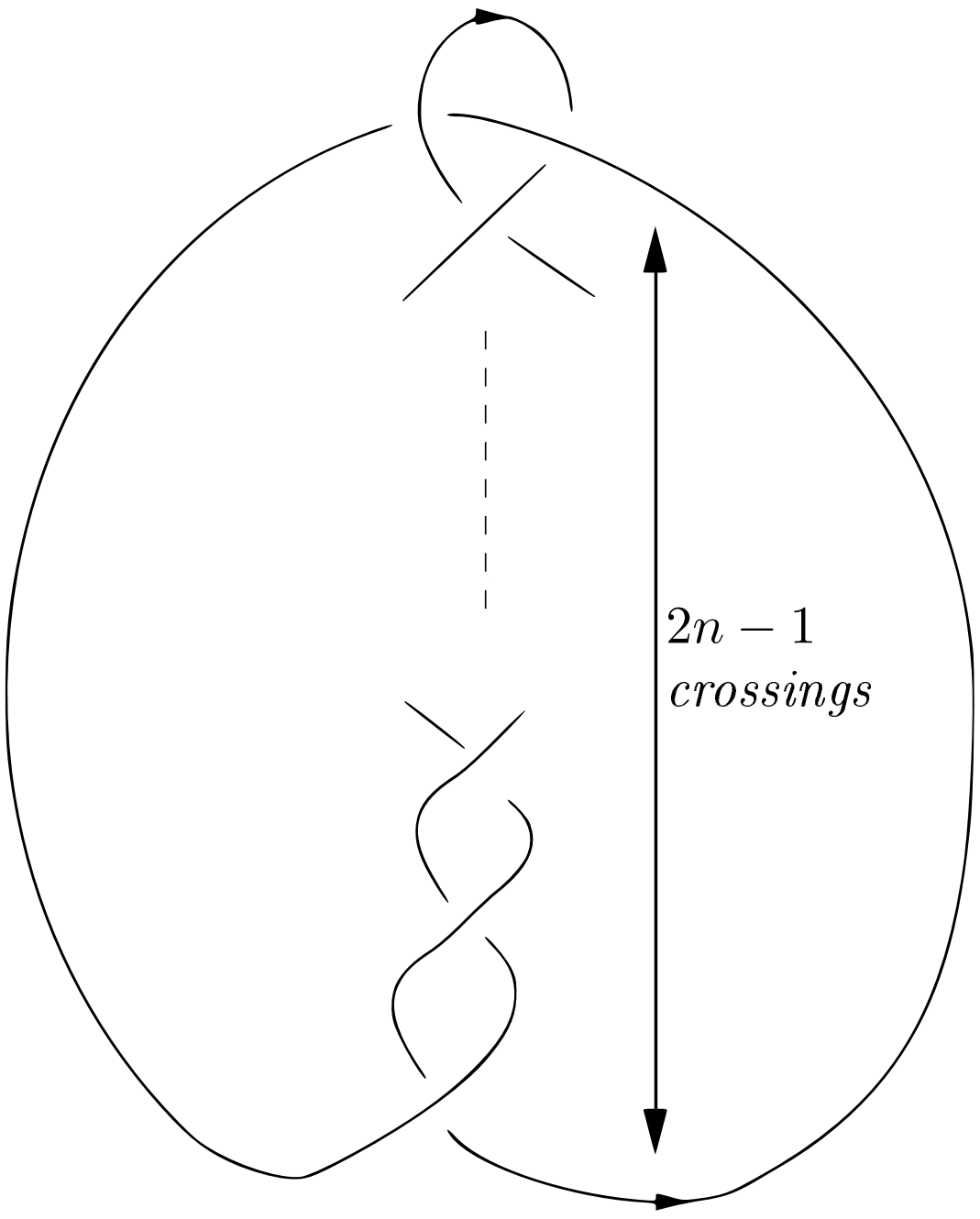

FIGURE 5.

Say $m$ is even.

$$\begin{aligned}LG(D(b_1,\dots,b_m)) =&a_1(b_m)LG(D(b_1,\dots,b_{m-1}-1))\\&+a_2(b_m)LG(D(b_1,\dots,b_{m-1}))\\&+a_3(b_m)LG(D(b_1,\dots,b_{m-2})).\end{aligned}$$

Since $LG(D(b_1,\dots,b_{m-1},0)) = LG(D(b_1,\dots,b_{m-2}))$ the first element in the sum has a span smaller than $0+2(m-1)=2m-2$. The second part has span $2+2(m-1)=2m$ and the third $2+2(m-2)=2m-2$. In the end

$$span(LG(D(b_1,\dots,b_m))) = 2m.$$

CQFD

2.3.3. *Twist knots.*

**Definition 2.16.** A twist knot is a Whitehead double of the unknot. We will denote by $K_n$ the twist knot shown in Figure 5 when $2n-1$ is positive. If $2n-1$ is negative, there are $1-2n$ crossings of the opposite sort.

For instance $K_0$ is the unknot, $K_1$ is the trefoil knot and $K_2$ is knot $5_2$.

**Proposition 2.17.** *For $n\neq 0$, $g(K_n)=1$.*

**Proposition 2.18.**

- $LG(K_0) = LG(\bigcirc) = 1.$
- *For $n\geqslant 1$,*

$$\begin{aligned}LG(K_n) =&\\&(-t_0^{-2}t_1^{-1}-t_0^{-1}t_1^{-2}+t_0^{-2}+2t_0^{-1}t_1^{-1}+t_1^{-2}-t_0^{-1}-t_1^{-1}+1)(t_0^{-1}t_1^{-1}\tilde{a_1}(n-1)+1)\\&+(t_0^{-1}-1)^2(t_1^{-1}-1)^2\frac{(t_0^{-1}t_1^{-1}+1)\tilde{a_1}(n-1)-2(n-1)}{t_0^{-1}t_1^{-1}-1}\\&+\tilde{a_1}(n-1)(t_0^{-1}-1)^2(t_1^{-1}-1)^2-t_0^{-1}t_1^{-1}\tilde{a_1}(n-1).\end{aligned}$$

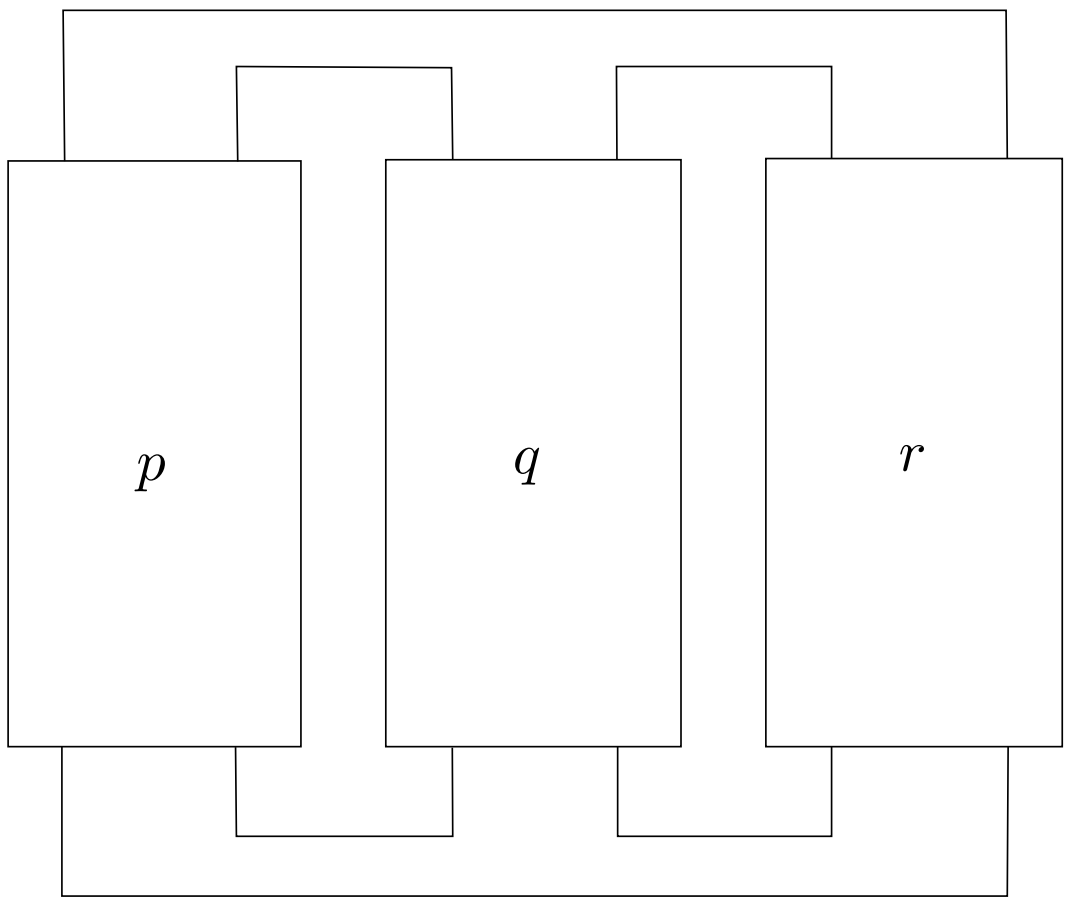

FIGURE 6.

*It has span 4.*

- *For* $n \leqslant -1$,

$$LG(K_n) = (t_0^{-1}-1)(t_1^{-1}-1)\left(a_1(-n)+(t_0-1)(t_1-1)\frac{2n+a_1(-n)(t_0t_1+1)}{t_0t_1-1}\right) + (t_0-1)(t_1-1)a_1(-n)+1.$$

*It has span 4 as well.*

*Proof.* When $n \geqslant 1$, we first write the third skein relation of 2.12 for $n-1$ full twists. On each of the three links that appear, we use the first point of 2.12. We find the formula written in the theorem. A close look at that expression shows that

$$span(LG(K_n)) \leqslant 4.$$

To see it is equal to 4, we can for example evaluate $t_1 = -t_0^{-1}$. We know we will find (and can verify)

$$LG(K_n)(t_0, -t_0^{-1}) = \Delta_{K_n}(t_0^2) = nt_0^2 - (2n-1) + nt_0^{-2}.$$

So $span(LG(K_n)) \geqslant 4$. Similar computations can be made when $n \leqslant -1$. CQFD

2.3.4. *Pretzel knots.*

**Definition 2.19.** Set $p, q, r \in \mathbb{Z}$. The $(p,q,r)$-pretzel link $L(p,q,r)$ is a union of three pairs of strands half-twisted $p, q, r$ times and attached along the tops and bottoms as shown in Figure 6. The half-twists are oriented according to whether the integer is positive or negative.

For example, pretzel knot $L(-2,3,7)$ is represented in Figure 7.

**Proposition 2.20.**

$$L(p,q,r) \text{ is a knot} \iff \big(\text{at most one of the three integers } p, q \text{ and } r \text{ is even}\big).$$

*In that case pretzel knot* $L(p,q,r)$ *is denoted by* $K(p,q,r)$.

In [15], Kim and Lee explicit the genus for all pretzel knots. Verifying the genus conjecture on this family of knots is quite interesting since the genus does not behave the same way as a function of parameters $(p,q,r)$ in all cases. The next theorem is proved in [15].

**Theorem 2.21.** *Let* $p, q, r$ *be integers. The genus of* $K(p,q,r)$ *is as follows :*

**1 :** $K(p,\pm1,\mp1)$, $K(\pm2,\mp1,\pm3)$ *have genus 0 for any* $p$,

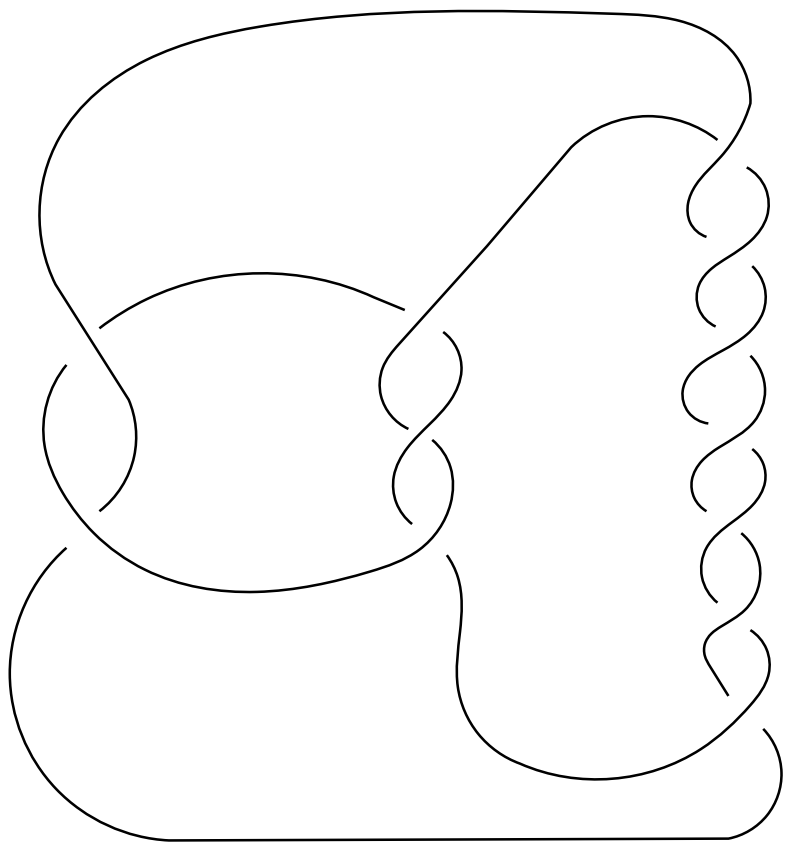

Figure 7. $L(-2,3,7)$.

**2 :** *$K(p,q,r)$ has genus 1 is $p,q,r$ are odd and we are not in case* **1***,*
**3 :** *$K(\pm 2,\mp 1,\pm r)$ has genus $\frac{|r-2|-1}{2}$,*
**4 :** *$K(2l,q,r)$ has genus $\frac{|q|+|r|}{2}$ if $q,r$ have the same sign and we are not in any of the previous cases,*
**5 :** *$K(2l,q,r)$ has genus $\frac{|q|+|r|-2}{2}$ if $q,r$ have different signs and we are not in cases* **1***,* **2** *or* **3***.*

We rewrite that theorem so that different cases exclude each other. Doing this makes computations more specific and somewhat easier in each case. Moreover, since $K(p,q,r)^* = K(-p,-q,-r)$, we will consider $p \geqslant 0$. Also, $K(p,q,r) = K(q,r,p) = K(r,p,q)$. So we can restrict our study to the cases where $q,r$ are odd.

**Corollary 2.22.** *Given the restrictions mentioned, setting $p \geqslant 0$ an integer and $q,r$ two odd integers, the genus $g$ of $K(p,q,r)$ is as follows :*

**1 :** *$g=0$ for $K(p,\pm 1,\mp 1)$ and $K(2,-1,3)$, that is when $K(p,q,r)$ is the unknot,*
**2 :** *$g=1$ when $p,q,r$ are odd and $K(p,q,r)$ is not the unknot,*
**3 :** *$g=\frac{|r-2|-1}{2}$ for $K(2,-1,r)$,*
**4 :** *$g=\frac{|q|+|r|}{2}$ if :*
- *$p$ is even and $q,r$ are positive,*
- *$p$ is even and different from 2, $q=-1$ and $r$ is negative,*
- *$p$ is even, $q$ is negative and different from $-1$, $r$ is negative,*

**5 :** *$g=\frac{|q|+|r|-2}{2}$ if :*
- *$p$ is even and different from 2, $q=-1$ and $r \geqslant 3$,*
- *$p$ is even, $q \geqslant 0$, $r \leqslant 0$ and $(p,q,r) \neq (p,1,-1)$,*
- *$p$ is even, $q \leqslant -3$ and $r \geqslant 0$.*

**Theorem 2.23.** *For all pretzel knots,*

$$span(LG(K(p,q,r))) \leqslant 4g(K(p,q,r)).$$

*Proof.* We compute the span of $LG(K(p,q,r))$ in each case of Corollary 2.22.
$\boxed{1}$ $K(p,\pm 1,\mp 1)$ and $K(2,-1,3)$ are different representations of the trivial knot that is part of the small cases we already checked.
$\boxed{2}$ Using the fact that $K(p,q,r)^* = K(-p,-q,-r)$ and $K(p,q,r) = K(q,r,p) = K(r,p,q)$, we have only two cases to consider : when $p,q,r$ have the same sign and when two out of the three have the same sign. For example we can choose the following configurations : $p,q,r \geqslant 0$ and $p \geqslant 0, q,r \leqslant 0$. In each case, using skein relation 2 of Corollary 2.12 on the three pairs of strands, we find a sum of 27 terms, each of which is symmetric of span

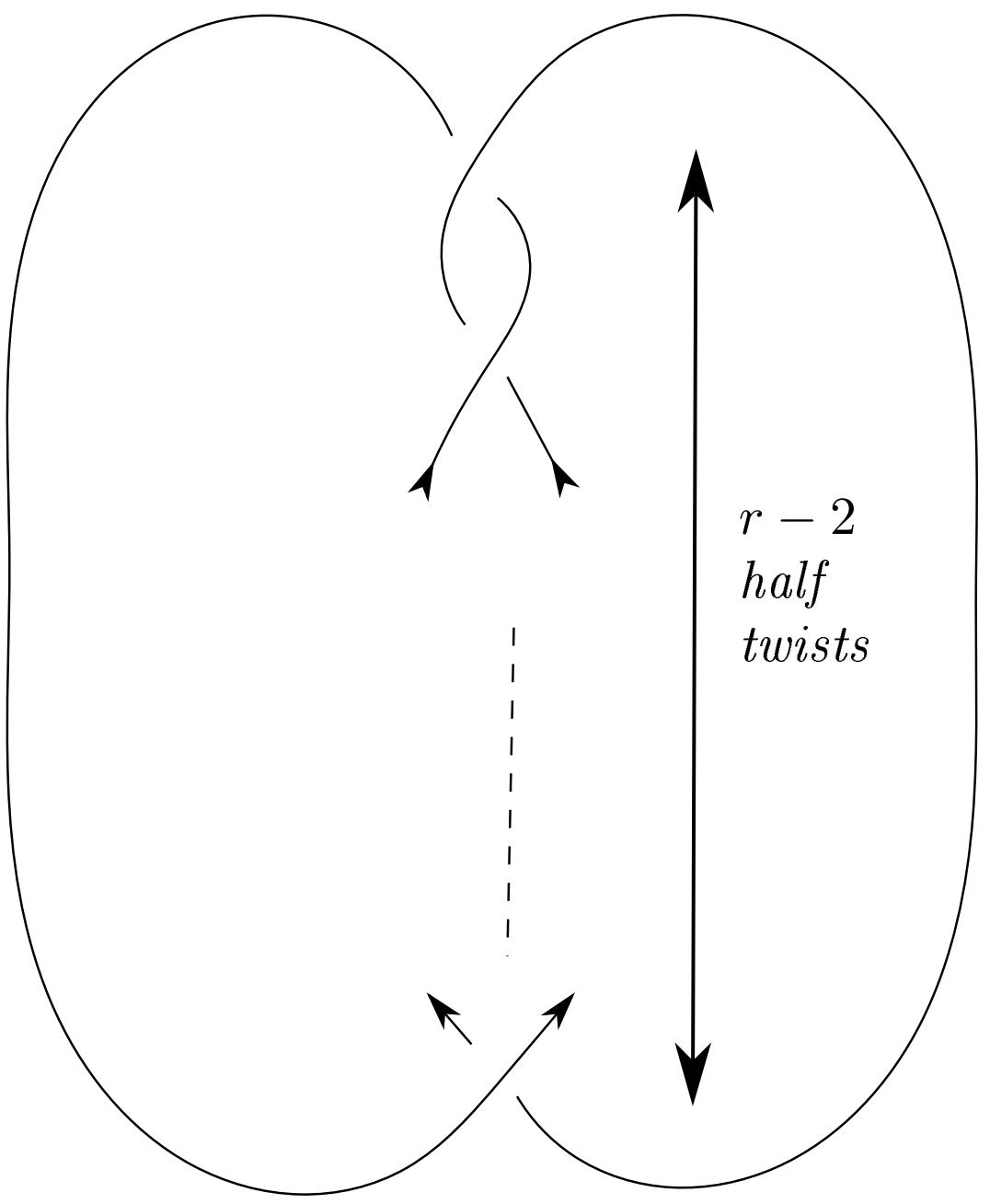

FIGURE 8. $K(2,-1,r)$ when $r \geqslant 3$.

smaller than 4.

**3** For $r = 1, 3$, $K(2,-1,r)$ is the unknot. If $r \geqslant 5$, $K(2,-1,r)$ is drawn in Figure 8 once it is simplified. The same kind of isotopy can be operated on $K(2,-1,r)$ when $r \leqslant -1$ and the result is shown in Figure 9.

For example if $r \geqslant 5$ we can apply skein relation 1 of Corollary 2.12 to the $r-2$ half twists.

$$\begin{aligned}
LG(K(2,-1,r)) &= \left(\frac{-1}{(t_0+1)(t_1+1)} + \frac{t_0^{r-2}}{(t_0+1)(t_0-t_1)} + \frac{t_1^{r-2}}{(t_1+1)(t_1-t_0)}\right) LG\left(\bigcirc\!\!\!\bigcirc\right) \\
&- \left(\frac{-t_0-t_1}{(t_0+1)(t_1+1)} + \frac{t_0^{r-2}(t_1-1)}{(t_0+1)(t_0-t_1)} + \frac{t_1^{r-2}(t_0-1)}{(t_1+1)(t_1-t_0)}\right) LG(\bigcirc) \\
&+ \left(\frac{-t_0t_1}{(t_0+1)(t_1+1)} - \frac{t_0^{r-2}t_1}{(t_0+1)(t_0-t_1)} - \frac{t_1^{r-2}t_0}{(t_1+1)(t_1-t_0)}\right) LG(\bigcirc\bigcirc) \\
&= \frac{(t_0-t_1)(t_0t_1+1) - t_0^{r-1}(t_1-1)(t_1+1) + t_1^{r-1}(t_0-1)(t_0+1)}{(t_0+1)(t_1+1)(t_0-t_1)}
\end{aligned}$$

The numerator has span $2r-2$ and the denominator has span 4. So

$$span(LG(K(2,-1,r))) = (2r-2) - 4 = 2r-6 = 4\left(\frac{r-3}{2}\right) = 4\left(\frac{|r-2|-1}{2}\right).$$

Computations can be led in a similar way in the other case.

**4** *If $p$ is positive and even, and $q, r$ are positive and odd.* Choosing an orientation for $K(p,q,r)$, we can apply skein relation 2 of 2.12 on the $p$ half twists and skein relation 1 of 2.12 on the $q$ and $r$ half twists. As in a previous case, we get a sum of 27 terms, each of which can be computed easily. All these terms have a span smaller than $2q + 2r = 4g$. Therefore we have the inequality in this case.

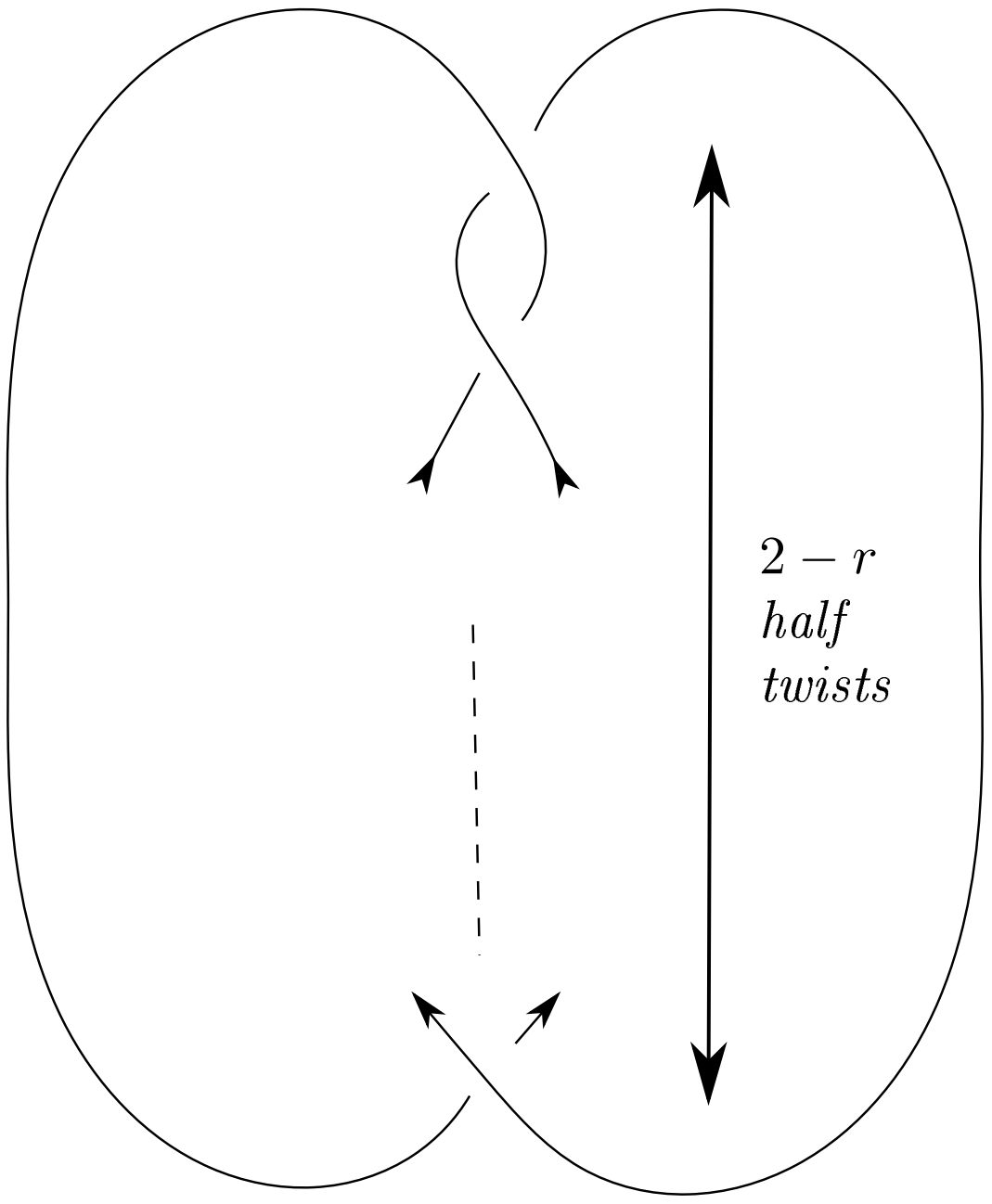

FIGURE 9. $K(2,-1,r)$ when $r \leqslant -1$.

*If $p$ is positive, even and different from 2, $q=-1$, and $r$ is odd and negative.* Choosing an orientation here again, we can use skein relation 2 on the $p$ half twists and skein relation 1 on the $r$ half twists. Each of the 9 parts of the sum such obtained has a span smaller than $2-2r$. So once again

$$span(LG(K(p,-1,r))) \leqslant 2-2r = 4\left(\frac{1-r}{2}\right) = 4g.$$

*If $p$ is even and positive, $q$ is odd negative and different from $-1$, and $r$ is odd and negative.* An extended computation similar to the two previous ones proves the bound in this case as well.

**5** *If $p$ is even, positive and different from 2, $q=-1$, and $r$ is positive, odd and different from* 1. Applying skein 2 of 2.12 on the $p/2$ full twists, we find three links, each of which is the closure of a power of generator $\sigma_1$ of $B_2$. We can then use skein relation 1 of 2.12 on each of these links to find that these three terms have a span smaller than $2r-2 = 4\left(\frac{r-1}{2}\right)$.

*If $p$ is even and positive, $q$ is odd and positive, $r$ is odd and negative, and $(q,r) \neq (1,-1)$.* This is the most tricky case. Indeed, if we compute $LG$ naively using skein relations 1 and 2 as we did for the moment, some parts of the sum we obtain have a span larger than $4g = 2q-2r-4$. We therefore have to look at these particular terms to see that what goes past the bound we hope actually compensates. This is achieved in the appendix (section 4).

*If $p$ is positive and even, $q$ is odd and $q \leqslant -3$, and $r$ is odd and positive.* $K(p,q,r)$ can be isotoped as follows :

$$K(p,q,r) = K(r,q,p) = K(p,r,q).$$

The last form $K(p,r,q)$ shows that this case is a consequence of the two previous cases of point **5**. CQFD

2.4. **A generalization of Conjecture 0.3 to a family of Links-Gould invariants.** The Links-Gould polynomial we defined and studied up to now is a particular case of a larger family of Links-Gould invariants, introduced by David de Wit in [7]. We will write $LG^{m,n}$ where $m, n$ are positive integers. Each invariant is derived from a highest weight representation of $U_q(gl(m|n))$. The invariant we denoted $LG$ is $LG^{2,1}$.

In [4], D. De Wit, A. Ishii and J. Links proved that

$$LG^{n,1}(L; t_0, e^{-2i\pi/n} t_0^{-1}) = \Delta_L(t_0^n) \text{ for any } n \in \mathbb{N}^*.$$

They also conjectured that for any $n \in \mathbb{N}^*$

$$LG^{n,1}(L; t_0, t_0^{-1}) = \Delta_L(t_0)^n.$$

This was proved by the author when $n = 2, 3$ by studying the link between the R-matrix representations of braid groups that give birth to the Links-Gould invariants and the Burau representations of $B_n$ [16]. More recently, Bertrand Patureau-Mirand and the author proved the statement for any $n$ by showing the $-1$ specialization of quantum supergroup $U_q(gl(n|1))$ shares properties with super Hopf algebra $U_q(gl(1|1))$ [17]. Though we lack computations for $LG^{m,n}$ when $(m, n) \neq (1, 1), (2, 1)$, the previous evaluations extend those we have for $LG^{2,1}$ so a potential homological interpretation for $LG^{2,1}$ should extend to $LG^{n,1}$ as well.

**Question 2.24.** *Set $L$ a link in $S^3$ and $n \geqslant 3$. Do we have, as it seems to be the case when $n = 2$ :*

- **I:** *$span(LG^{n,1}(L; t_0, t_1)) \leqslant n(2g(L) + \mu - 1)$,*
- **II:** *If $L$ is alternating, then inequality* **I** *is an equality ?*

*Remark* 2.25. For example, the equality holds for all prime knots with less than 10 crossings when $n = 3$.

*Remark* 2.26. As a consequence, one could ask, as $n$ tends to infinity, if

$$span(LG^{n,1}(L; t_0, t_1)) \underset{+\infty}{\sim} n(2g(L) + \mu - 1) \text{ ?}$$

However this *cannot* be true. Indeed, there are pairs of mutant knots with different genera, and neither of the $LG^{n,1}$ detects mutation.

## 3. The Links-Gould polynomial and fiberedness

In addition to genus information, the Links-Gould polynomial seems to contain signs of whether a knot is fibered or not. This is another well known feature of the Alexander invariant. Conjecture 0.4, if it were to be true, would refine the standard Alexander polynomial criterion. This is the object of this section.

**Definition 3.1.** A knot $K$ in $S^3$ is said to be fibered if the two following conditions hold :

- **1:** the complement of the knot is the total space of a locally trivial bundle over the base space $S^1$, i.e. there exists a map $p : S^3 \setminus K \longrightarrow S^1$ which is a locally trivial bundle.
- **2:** there exists $V(K)$ a neighborhood of $K$ and there exists a trivializing homeomorphism $\theta : V(K) \longrightarrow S^1 \times D^2$ such that $\pi \circ \theta(X) = p(X)$ for any $X \in V(K) \setminus K$, where $\pi(x, y) := \frac{y}{|y|}$.

Let us recall well-known properties of the Alexander polynomial of a fibered knot. The Alexander polynomial of a fibered knot is monic [22, 24, 28]. This means the coefficient of the highest degree term of the standard Alexander normalization of the polynomial is 1. For the Conway normalization, it means the leading coefficient is $\pm 1$. The converse is not true in general. However, the condition is sufficient for prime knots with up to 10 crossings and alternating knots [21]. Also note that for fibered knots, the degree of the Alexander polynomial is exactly twice the genus of the knot, that is the genus of the corresponding fibre surface [24]. This last property yields the following.

**Proposition 3.2.** *Set $K$ a fibered knot. If Conjecture 0.3 is true, then*

$$span(LG(K;t_0,t_1)) = 4g(K).$$

*Proof.* To have such a result we can more generally consider a set of links $E$ such that, for any $L \in E$, $deg(\Delta_L(t)) = 2g(L) + \mu - 1$. This is the case here, and it was also the case in Proposition 1.13 where it is proved completely. CQFD

We can express the Links-Gould polynomial with different sets of variables : $(t_0, t_1)$ as we did up to now, but also $(p, q)$ where

$$p^a q^b = t_0^{-\frac{a}{4}+\frac{b}{2}} t_1^{\frac{a}{4}+\frac{b}{2}}.$$

These are the variables used in the LINKS-GOULD EXPLORER as well as in de Wit's papers on the subject. He sometimes uses $P = p^2$. In variables $(p, q)$ the $LG$ polynomial of a link $L$ can be written

$$LG(L;p,q) = a_0 + \sum_{k \in \mathbb{N}^*} P_k(q)(p^{2k} + p^{-2k})$$

where $a_0 \in \mathbb{Z}$ and $P_k(q) \in \mathbb{Z}[q^{\pm 1}]$. Note that if $P_l(q) \neq 0$ and $P_k(q) = 0$ for any $k > l$, then $span(LG(L;p,q)) = 2l$.

**Definition 3.3.** Set $K$ a knot. We say $LG(K)$ is *monic* when the term in $LG(K)$ of highest and lowest degrees can be written $q^{2m}(p^{4l} + p^{-4l})$ with $l \in \mathbb{N}$ and $m \in \mathbb{Z}$. In terms of variables $(t_0, t_1)$, this condition is expressed by saying the terms of highest and lowest degrees are monic monomials of the form $t_0^{\alpha} t_1^{\beta}$ with $\alpha + \beta$ even.

**Proposition 3.4.** *Set $K$ a knot. If $LG(K;t_0,t_1)$ is monic, then $\Delta_K(t)$ is monic as well.*

*Proof.* Consequence of $LG(K;t,-t^{-1}) = \Delta_K(t^2)$. CQFD

*Remark* 3.5. Point **I** in Conjecture 0.4 implies point **II**. This is a consequence of Proposition 3.4 and of the fact that when the Alexander polynomial is monic for an alternating knot, the knot is fibered.

*Remark* 3.6. Given Proposition 3.4, criterion 0.4 would be an improvement of the criterion provided by the Alexander invariant. In addition, we will see in the following that there are examples of knots where $\Delta$ is monic but $LG$ is not.

**Proposition 3.7.**

**1 :** *Conjecture 0.4 holds for every prime knot up to 12 crossings. In particular, for all alternating knots tested, fiberedness and having monic LG are equivalent.*

**2 :** *For a prime knot $K$ with at most 11 crossings, $K$ is fibered if and only if $LG(K)$ is monic.*

By work of Friedl and Kim [9], there are 13 non-fibered 12-crossing prime knots which have monic Alexander polynomials such that $deg(\Delta_K(t)) = 2g(K)$ : $12^N_{57}$, $12^N_{210}$, $12^N_{214}$, $12^N_{258}$, $12^N_{279}$, $12^N_{382}$, $12^N_{394}$, $12^N_{464}$, $12^N_{483}$, $12^N_{535}$, $12^N_{650}$, $12^N_{801}$, $12^N_{815}$. Among these, $LG$ manages to "detect" non-fiberedness of some, but not all.

**Proposition 3.8.** *The following knots have monic $LG$ : $12^N_{57}$, $12^N_{258}$, $12^N_{279}$, $12^N_{464}$, $12^N_{483}$, $12^N_{650}$, $12^N_{815}$. So there are knots that are non-fibered but that have monic $LG$.*

**Proposition 3.9.** *The following knots have non-monic $LG$ : $12^N_{210}$, $12^N_{214}$, $12^N_{382}$, $12^N_{394}$, $12^N_{535}$, $12^N_{801}$. So $\Delta$ sometimes is monic when $LG$ is not.*

Verifications were done using de Wit's LINKS-GOULD EXPLORER [8] and Cha and Livingston's KNOTINFO [2].

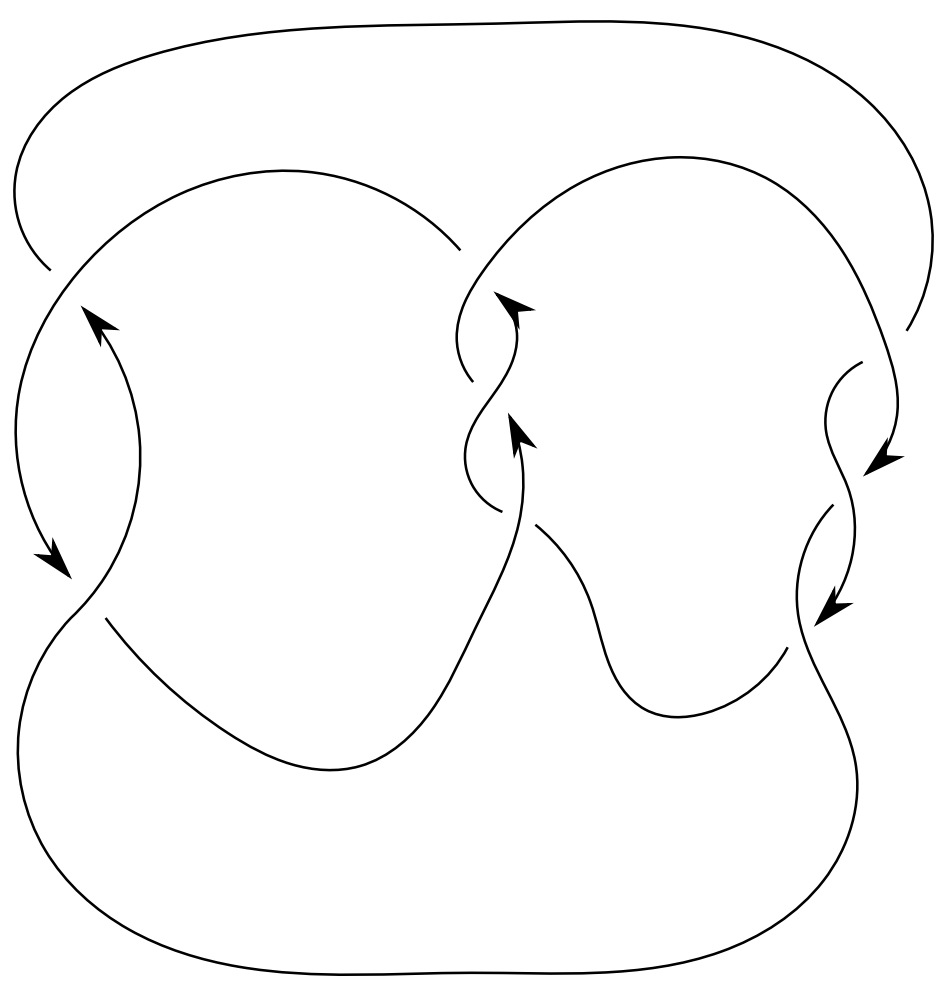

FIGURE 10. $K(2,3,-3)$.

## 4. Appendix : proof of the harder case in Theorem 2.23

Here we prove the remaining case of Theorem 2.23. Consider pretzel knot $K(p,q,r)$ when $p$ is positive and even, $q$ is positive and odd, $r$ is negative and odd, and $(q,r) \neq (1,-1)$. We want to show that

$$span(LG(K(p,q,r))) \leqslant 4g = 2q - 2r - 4.$$

We first consider $r=-1$. Then $4g = 2q-2$. In that precise configuration, using skein relation 2 of 2.12 on the $p/2$ full twists and skein relation 1 of the same corollary on the $q$ half twists, you find a sum of terms, each of which has a span smaller than $2q-2$.

In general, the computation is not that easy. We show knot $K(2,3,-3)$ in Figure 10 to fix the orientation chosen here. Let us introduce some notations :

$$x(n) = \frac{(-1)^n}{(t_0+1)(t_1+1)} + \frac{t_0^n}{(t_0+1)(t_0-t_1)} + \frac{t_1^n}{(t_1+1)(t_1-t_0)},$$

$$y(n) = \frac{(-1)^n(t_0+t_1)}{(t_0+1)(t_1+1)} + \frac{t_0^n(t_1-1)}{(t_0+1)(t_0-t_1)} + \frac{t_1^n(t_0-1)}{(t_1+1)(t_1-t_0)},$$

$$z(n) = \frac{(-1)^n t_0 t_1}{(t_0+1)(t_1+1)} - \frac{t_0^n t_1}{(t_0+1)(t_0-t_1)} - \frac{t_1^n t_0}{(t_1+1)(t_1-t_0)}.$$

We transform the $-r$ half twists in $K(p,q,r)$ using skein relation 1 of 2.12.

$$\begin{aligned} LG(K(p,q,r)) &= x(-r)(t_0^{-1},t_1^{-1})LG(K(p,q,-2)) \\ &\quad - y(-r)(t_0^{-1},t_1^{-1})LG(K(p,q,-1)) \\ &\quad + z(-r)(t_0^{-1},t_1^{-1})LG(K(p,q,0)). \end{aligned}$$

The span of $LG(K(p,q,-1))$ is $2q-2$ with the previous case. Also $span(y(-r)(t_0^{-1},t_1^{-1})) = -2r-4$. So the second term of the sum above has span $2q-2r-6$ and we need only to consider the first and third terms in the rest of the proof.

Also, using skein relation 1 of Proposition 2.7,

$$\begin{aligned} x(-r)(t_0^{-1},t_1^{-1})LG(K(p,q,-2)) &= x(-r)(t_0^{-1},t_1^{-1})(t_0^{-1}+t_1^{-1}-1)LG(K(p,q,-1)) \\ &\quad + x(-r)(t_0^{-1},t_1^{-1})(t_0^{-1}+t_1^{-1}-t_0^{-1}t_1^{-1})LG(K(p,q,0)) \\ &\quad - x(-r)(t_0^{-1},t_1^{-1})t_0^{-1}t_1^{-1}LG(K(p,q,1)). \end{aligned}$$

Again, the first term of this sum has span $(-2r-4)+2+(2q-2) = 2q-2r-4$ so we can ignore it as well and we are interested in the following sum :

$$\Big(x(-r)(t_0^{-1},t_1^{-1})(t_0^{-1}+t_1^{-1}-t_0^{-1}t_1^{-1}) + z(-r)(t_0^{-1},t_1^{-1})\Big)LG(K(p,q,0)) \\ - \Big(x(-r)(t_0^{-1},t_1^{-1})t_0^{-1}t_1^{-1}\Big)LG(K(p,q,1)).$$

However, $span(z(-r)(t_0^{-1},t_1^{-1})) = -2r-6$ and $span(x(-r)(t_0^{-1},t_1^{-1})(t_0^{-1}+t_1^{-1}-t_0^{-1}t_1^{-1})) = (-2r-4)+2 = -2r-2$. So we can reduce our concerns to

$$\Big(x(-r)(t_0^{-1},t_1^{-1})(t_0^{-1}+t_1^{-1}-t_0^{-1}t_1^{-1})\Big)LG(K(p,q,0)) \\ - \Big(x(-r)(t_0^{-1},t_1^{-1})t_0^{-1}t_1^{-1}\Big)LG(K(p,q,1)).$$

Using the usual skein relations on $LG(K(p,q,0)$, we get the following value modulo terms with a small span :

$$LG(K(p,q,0) = \Big(\frac{(t_0^{-1}t_1^{-1})^{p/2}-1}{t_0^{-1}t_1^{-1}-1}\big(-(t_0^{-1}-1)(t_1^{-1}-1)\big)x(q)(t_0+t_1-1-t_0t_1)\Big) \\ + \Big(\frac{2(t_0^{-1}-1)(t_1^{-1}-1)}{t_0^{-1}t_1^{-1}-1}\Big)\Big(\frac{p}{2}-\frac{(t_0^{-1}t_1^{-1})^{p/2}-1}{t_0^{-1}t_1^{-1}-1}\Big)x(q)\big(-(t_0-1)(t_1-1)\big).$$

The two pieces of the sum have span $2q$. Similarily, ignoring non extremal span terms :

$$LG(K(p,q,1) = \\ \Big(\frac{(t_0^{-1}t_1^{-1})^{p/2}-1}{t_0^{-1}t_1^{-1}-1}x(q)(t_0+t_1-1)(t_0^{-1}+t_1^{-1}-2-t_0^{-1}t_1^{-1})\big(-(t_0-1)(t_1-1)\big)\Big) \\ + \Big(\frac{2(t_0^{-1}-1)(t_1^{-1}-1)}{t_0^{-1}t_1^{-1}-1}\Big)\Big(\frac{p}{2}-\frac{(t_0^{-1}t_1^{-1})^{p/2}-1}{t_0^{-1}t_1^{-1}-1}\Big)x(q+1)\big(-(t_0-1)(t_1-1)\big).$$

So the quantity we are interested in is

$$
\begin{aligned}
x(-r)(t_0^{-1},t_1^{-1})(t_0^{-1}+t_1^{-1}-t_0^{-1}t_1^{-1})\Bigg[&\frac{(t_0^{-1}t_1^{-1})^{p/2}-1}{t_0^{-1}t_1^{-1}-1}\big(-(t_0^{-1}-1)(t_1^{-1}-1)\big)x(q)(t_0+t_1)\\
&+\frac{2(t_0^{-1}-1)(t_1^{-1}-1)}{t_0^{-1}t_1^{-1}-1}\left(\frac{p}{2}-\frac{(t_0^{-1}t_1^{-1})^{p/2}-1}{t_0^{-1}t_1^{-1}-1}\right)x(q)(t_0+t_1)\Bigg]\\
&+x(-r)(t_0^{-1},t_1^{-1})(-t_0^{-1}t_1^{-1})\Bigg[\frac{(t_0^{-1}t_1^{-1})^{p/2}-1}{t_0^{-1}t_1^{-1}-1}x(q)(t_0+t_1)(t_0^{-1}+t_1^{-1})(t_0+t_1)\\
&+\frac{2(t_0^{-1}-1)(t_1^{-1}-1)}{t_0^{-1}t_1^{-1}-1}\left(\frac{p}{2}-\frac{(t_0^{-1}t_1^{-1})^{p/2}-1}{t_0^{-1}t_1^{-1}-1}\right)x(q+1)(t_0+t_1)\Bigg]\\
=x(-r)(t_0^{-1},{t_1}^{-1})(t_0+t_1)\Bigg[&(t_0^{-1}+t_1^{-1})^2\frac{(t_0^{-1}t_1^{-1})^{p/2}-1}{t_0^{-1}t_1^{-1}-1}x(q)\\
&+(t_0^{-1}+t_1^{-1})\frac{2(t_0^{-1}-1)(t_1^{-1}-1)}{t_0^{-1}t_1^{-1}-1}\left(\frac{p}{2}-\frac{(t_0^{-1}t_1^{-1})^{p/2}-1}{t_0^{-1}t_1^{-1}-1}\right)x(q)\\
&+\frac{(t_0^{-1}t_1^{-1})^{p/2}-1}{t_0^{-1}t_1^{-1}-1}x(q)(t_0+t_1)(t_0^{-1}+t_1^{-1})(-t_0^{-1}t_1^{-1})\\
&+(-t_0^{-1}t_1^{-1})\frac{2(t_0^{-1}-1)(t_1^{-1}-1)}{t_0^{-1}t_1^{-1}-1}\left(\frac{p}{2}-\frac{(t_0^{-1}t_1^{-1})^{p/2}-1}{t_0^{-1}t_1^{-1}-1}\right)x(q+1)\Bigg]\\
=x(-r)(t_0^{-1},{t_1}^{-1})(t_0+t_1)\Bigg[&x(q)\frac{(t_0^{-1}t_1^{-1})^{p/2}-1}{t_0^{-1}t_1^{-1}-1}\big(t_0^{-2}+t_1^{-2}-t_0^{-1}t_1^{-1}(t_0t_1^{-1}+t_1t_0^{-1})\big)\\
&+\frac{2(t_0^{-1}-1)(t_1^{-1}-1)}{t_0^{-1}t_1^{-1}-1}\left(\frac{p}{2}-\frac{(t_0^{-1}t_1^{-1})^{p/2}-1}{t_0^{-1}t_1^{-1}-1}\right)\big((t_0^{-1}+t_1^{-1})x(q)-t_0^{-1}t_1^{-1}x(q+1)\big)\Bigg].
\end{aligned}
$$

But $t_0^{-2}+t_1^{-2}-t_0^{-1}t_1^{-1}(t_0t_1^{-1}+t_1t_0^{-1})=0$. So to show the two terms of highest and lowest degree disappear in that polynomial we show that modulo lower degree terms, $\alpha=(t_0^{-1}+t_1^{-1})x(q)-t_0^{-1}t_1^{-1}x(q+1)=0$. Let's look at $x(q)$ first of all :

$$x(q)=\frac{t_0^q-t_1^q+\text{ other terms}}{(t_0+1)(t_1+1)(t_0-t_1)}=M+m+\text{ other terms},$$

where $M$ is the term of highest degree in $x(q)$, and $m$ the one of smallest degree. That way

$$t_0^q-t_1^q+\text{ other terms }=(t_0+1)(t_1+1)(t_0-t_1)(M+m+\text{ other terms}).$$

And identifying the highest and lowest degree terms on each side we find

$$M=t_0^{q-2}\text{ and }m=t_1^{q-2}.$$

Finally, modulo lower degree terms,

$$
\begin{aligned}
\alpha&=(t_0^{-1}+t_1^{-1})(t_0^{q-2}+t_1^{q-2})-t_0^{-1}t_1^{-1}(t_0^{q-1}+t_1^{q-1})\\
&=t_1^{-1}t_0^{q-2}-t_0^{-1}t_1^{-1}t_0^{q-1}+t_0^{-1}t_1^{q-2}-t_0^{-1}t_1^{-1}t_1^{q-1}\\
&=0.
\end{aligned}
$$

In conclusion,

$$span(LG(K(p,q,r)))\leqslant(2q-2r-2)-\underline{2}=2q-2r-4=4\left(\frac{q-r-2}{2}\right).$$

Acknowledgments. I owe many warm thanks to my advisor Emmanuel Wagner who taught me about Links-Gould invariants and has been of much help in designing and writing this paper. I am most grateful to Hugh Morton for his help during my study of the Whitehead doubles of the trefoil, and for showing me a braid presentation for the untwisted double. I am indebted to David de Wit for the extended work he did implementing the LINKS-GOULD EXPLORER, and to Jon Links for telling me about this remarkable package. Finally I would like to thank David Cimasoni and Anthony Conway for the fruitful and pleasant conversations, as well as for their comments and ideas.

## References

[1] J. W. Alexander. *Topological invariants of knots and links.* Trans. Amer. Math. Soc., 30 (2): 275–306, 1928.

[2] J. C. Cha, C. Livingston. *KnotInfo: Table of Knot Invariants.* http://www.indiana.edu/~knotinfo

[3] R. H. Crowell. *Genus of Alternating Link Types.* Ann. of Math., (2) Vol. 69 (1959), 258-275.

[4] D. De Wit, A. Ishii, J. Links. *Infinitely many two-variable generalisations of the Alexander-Conway polynomial.* Algebraic and Geometric Topology, Volume 5, 405-418, 2005.

[5] D. De Wit, L.H. Kauffman, J. Links. *On the Links-Gould invariant of links.* J. Knot Theory Ramifications. 8 (1999), no.2, 165-199. MR 2000j:57020

[6] D. De Wit, J. Links. *Where the Links-Gould invariant first fails to distinguish nonmutant prime knots.* J. Knot Theory Ramifications. 16, 1021 (2007). DOI: 10.1142/S0218216507005658.

[7] D. De Wit. *An infinite suite of Links-Gould invariants.* J. Knot Theory Ramifications. 10 (2001), no.1, 37-62. MR 2002e:57015

[8] D. De Wit. *Links-Gould Explorer.* http://www.maths.uq.edu.au/cmp/Links--Gould%20Explorer.html

[9] S. Friedl, T. Kim. *The Thurston norm, fibered manifolds and twisted Alexander polynomials.* Topology, 45, 929-953 (2006).

[10] N. Geer, B. Patureau-Mirand. *Multivariable link invariants arising from sl(2|1) and the Alexander polynomial.* J. Pure Appl. Algebr., Volume 210, Issue 1, July 2007, Pages 283–298.

[11] J. Hoste, M. Thistlethwaite, J. Weeks. *The first 1,701,936 knots.* The Mathematical Intelligencer. 20(4):33–48, 1998.

[12] A. Ishii. *The Links-Gould polynomial as a generalization of the Alexander-Conway polynomial.* Pacific J. Math., 225 (2006) 273–285.

[13] A. Ishii. *Algebraic links and skein relations of the Links-Gould invariant.* Proc. Amer. Math. Soc., 132 (2004), pp. 3741-3749.

[14] T. Kanenobu. *Genus and Kauffman polynomial of a 2-bridge knot.* Osaka J. Math., Volume 29, Number 3 (1992), 635-651.

[15] D. Kim, J. Lee. *Some invariants of pretzel links.* Bull. Austral. Math. Soc., Vol. 75 (2007), 253-271.

[16] B.-M. Kohli. *On the Links-Gould invariant and the square of the Alexander polynomial.* J. Knot Theory Ramifications, Vol. 25, No. 02, 1650006 (2016). DOI: 10.1142/S0218216516500061.

[17] B.-M. Kohli, B. Patureau-Mirand. *Other quantum relatives of the Alexander polynomial through the Links-Gould invariant.* Work in progress.

[18] J. Links, M.D. Gould. *Two variable link polynomials from quantum supergroups.* Letters in Mathematical Physics, 26(3):187-198, November 1992.

[19] H. R. Morton. *Seifert circles and knot polynomials.* Proc. Camb. Phil. Soc., 99 (1986), 107-109.

[20] K. Murasugi. *On the genus of the alternating knot II.* J. Math. Soc. Japan, Volume 10, Number 3 (1958), 235-248.

[21] K. Murasugi. *On a certain subgroup of the group of an alternating link.* Amer. J. Math., 85 (1963), 544-550.

[22] L. Neuwirth. *Knot Groups.* Ann. of Math. Stud. 56, Princeton University Press, Princeton, N.J., 1965.

[23] T. Ohtsuki. *Quantum invariants. A study of knots, 3-manifolds, and their sets.* Number 29 in Series on Knots and Everything. World Scientific, 2002.

[24] E. Rapaport. *On the commutator subgroup of a knot group.* Ann. of Math., (2) 71 (1960), 157-162.

[25] N. Reshetikhin, C. Stroppel, B. Webster. *Schur-Weyl-Type duality for quantized* $gl(1|1)$*, the Burau representation of braid groups, and invariants of tangled graphs.* I. Itenberg et al. (eds.), Perspectives

in Analysis, Geometry, and Topology: On the Occasion of the 60th Birthday of Oleg Viro, Progress in Mathematics 296.

[26] H. Seifert. *Über das Geschlecht von Knoten*. Math. Annalen, 110 (1): 571–592, 1934.

[27] L. Siebenmann. *Exercices sur les noeuds rationnels*. unpublished notes from the Département de Mathématiques d'Orsay, 1975.

[28] J. Stallings. *On fibering certain 3-manifolds*. In Topology of 3-manifolds and related topics (Proc. The Univ. of Georgia Institute, 1961), Prentice-Hall, Englewood Cliffs, N.J., 1962, 95-100.

[29] A. Stoimenow. *Knot Data Tables*. http://stoimenov.net/stoimeno/homepage/ptab/13n.gen

[30] O. Viro. *Quantum relatives of the Alexander polynomial*. St. Petersburg Math. J., Vol. 18 (2007), no.3, 391-457.

(Ben-Michael Kohli) IMB UMR5584, CNRS, Université Bourgogne Franche-Comté, F-21000 Dijon, France.

*E-mail address*: Ben-Michael.Kohli@u-bourgogne.fr